\documentclass[12pt]{extarticle}
\usepackage{mathptmx}
\usepackage{bbm}

\newcommand{\ind}{\mbox{\mbox{ } \mbox{ } \mbox{ }}}

\newcommand{\be}{\begin{enumerate}}
\newcommand{\ee}{\end{enumerate}}
\newcommand{\bc}{\begin{center}}
\newcommand{\ec}{\end{center}}
\newcommand{\beq}{\begin{equation}}
\newcommand{\eeq}{\end{equation}}
\newcommand{\bqn}{\begin{eqnarray}}
\newcommand{\eqn}{\end{eqnarray}}
\newcommand{\bqns}{\begin{eqnarray*}}
\newcommand{\eqns}{\end{eqnarray*}}
\newcommand{\N}{{\bf N}}
\newcommand{\R}{{\bf R}}

\newcommand{\rarrow}{\rightarrow}

\newcommand{\se}{{\cal E}}

\newcommand{\bdoc}{\begin{document}}
\newcommand{\edoc}{\end{document}}
\newcommand{\docart}{\documentstyle[12pt]{article}}
\newcommand{\vs}{\vspace{.16in}}
\newcommand{\half}{\frac{1}{2}}
\newcommand{\diffeo}{\mbox{ diffeomorphism }} 
\newcommand{\dnorm}[1]{\mbox{$\mid \mid #1 \mid \mid$}}
\newcommand{\twoover}[2]{{\footnotesize \begin{array}[t]{c}\mbox{ {\normalsize #1}}  \\
                               #2 
        \end{array} } \hspace{1cm} }
\newcommand{\threeover}[3]{ {\footnotesize \begin{array}[t]{c}\mbox{ {\normalsize #1}}  \\
                               #2  \\
                                 #3
        \end{array} } \hspace{1cm} }
\newcommand{\fourover}[4]{ {\footnotesize \begin{array}[t]{c}\mbox{ {\normalsize #1}}  \\
                               #2  \\
                                 #3 \\
				#4
        \end{array} } \hspace{1cm}  }

\newenvironment{vover3}[3]{   
       {\footnotesize \begin{array}[t]{c}
                               {\normalsize #1} \vspace{.2cm} \\ 
                                 #2 \\
                                 #3  
                       \end{array} }
                            } { \hspace{1cm} } 
\newenvironment{vover4}[4] {
       {\footnotesize \begin{array}[t]{c} 
                               {\normalsize #1} \vspace{.2cm} \\  
                                 #2 \\
                                 #3 \\
                                 #4  
                       \end{array} }
                            } {  \hspace{1cm} }

\newenvironment{seteqover}[3]{#1 \ \ = \hspace{.6cm}
                                #2 #3} {\ind} 

\newcommand{\norm}[1]{\mbox{$\mid #1 \mid$}}
\newcommand{\Norm}[1]{\mbox{$\parallel #1 \parallel$}}
\newcommand{\imp}{\mbox{$\Longrightarrow$}}
\newcommand{\pr}{\mbox{$\prime$}}
\newcommand{\br}{\mbox{\bf R}}
\newcommand{\fact}{\mbox{$\verb+!+$}}

\newcommand{\tag}[2]{\vspace{.15in}$(#1) \hfill #2 \hfill $ \vspace{.15in} }

\newcommand{\ba}{\mbox{${\bf a}$}}
\newcommand{\bb}{\mbox{${\bf b}$}}
\newcommand{\bi}{\mbox{${\bf i}$}}
\newcommand{\bj}{\mbox{${\bf j}$}}
\newcommand{\bfi}{\mbox{${\bf i}$}}
\newcommand{\bfj}{\mbox{${\bf j}$}}
\newcommand{\bfx}{\mbox{${\bf x}$}}
\newcommand{\bfy}{\mbox{${\bf y}$}}
\newcommand{\bfg}{\mbox{${\bf g}$}}
\newcommand{\bfc}{\mbox{${\bf c}$}}
\newcommand{\bS}{\mbox{${\bf S}$}}
\newcommand{\bV}{\mbox{${\bf V}$}}
\newcommand{\fh}{\hat{f}}
\newcommand{\Ih}{\hat{I}}
\newcommand{\Lh}{\hat{L}}

\newcommand{\alp}{\mbox{$\alpha$}}
\newcommand{\gra}{\mbox{$\alpha$}}

\newcommand{\bet}{\mbox{$\beta$}}
\newcommand{\grb}{\mbox{$\beta$}}

\newcommand{\gam}{\mbox{$\gamma$}}
\newcommand{\Gam}{\mbox{$\Gamma$}}
\newcommand{\del}{\mbox{$\delta$}}
\newcommand{\Del}{\mbox{$\Delta$}}
\newcommand{\sig}{\mbox{$\sigma$}}
\newcommand{\Sig}{\mbox{$\Sigma$}}

\newcommand{\grg}{\mbox{$\gamma$}}
\newcommand{\Grg}{\mbox{$\Gamma$}}
\newcommand{\grd}{\mbox{$\delta$}}
\newcommand{\Grd}{\mbox{$\Delta$}}
\newcommand{\grs}{\mbox{$\sigma$}}
\newcommand{\Grs}{\mbox{$\Sigma$}}
\newcommand{\grf}{\mbox{$\phi$}}
\newcommand{\Grf}{\mbox{$\Phi$}}
\newcommand{\grth}{\mbox{$\theta$}}
\newcommand{\Grth}{\mbox{$\Theta$}}
\newcommand{\gro}{\mbox{$\omega$}}
\newcommand{\Gro}{\mbox{$\Omega$}}
\newcommand{\gre}{\mbox{$\epsilon$}}
\newcommand{\grve}{\mbox{$\varepsilon$}}
\newcommand{\grr}{\mbox{$\rho$}}
\newcommand{\Grr}{\mbox{$\Rho$}}
\newcommand{\grt}{\mbox{$\tau$}}
\newcommand{\gri}{\mbox{$\iota$}}
\newcommand{\grp}{\mbox{$\pi$}}
\newcommand{\grP}{\mbox{$\Pi$}}
\newcommand{\grl}{\mbox{$\lambda$}}
\newcommand{\grL}{\mbox{$\Lambda$}}
\newcommand{\Grl}{\mbox{$\Lambda$}}
\newcommand{\grz}{\mbox{$\zeta$}}

\newcommand{\thet}{\mbox{$\theta$}}
\newcommand{\Thet}{\mbox{$\Theta$}}
\newcommand{\ome}{\mbox{$\omega$}}
\newcommand{\Ome}{\mbox{$\Omega$}}
\newcommand{\eps}{\mbox{$\epsilon$}}

\newcommand{\iot}{\mbox{$\iota$}}

\newcommand{\lam}{\mbox{$\lambda$}}
\newcommand{\Lam}{\mbox{$\Lambda$}}
\newcommand{\zet}{\mbox{$\zeta$}}

\newcommand{\kap}{\mbox{$\kappa$}}

\newcommand{\grx}{\mbox{$\xi$}}
\newcommand{\Grx}{\mbox{$\Xi$}}
\newcommand{\grc}{\mbox{$\chi$}}
\newcommand{\Grc}{\mbox{$\Chi$}}
\newcommand{\grn}{\mbox{$\nu$}}
\newcommand{\grm}{\mbox{$\mu$}}
\newcommand{\grk}{\mbox{$\kappa$}}

\newtheorem{theorem}{Theorem}[section]
\newtheorem{Cor}[theorem]{Corollary}
\newtheorem{Prop}[theorem]{Proposition}
\newtheorem{Lem}[theorem]{Lemma}
\newtheorem{Rem}[theorem]{Remark}
\newtheorem{Proof}[theorem]{Proof}

\newcommand{\rar}{\rightarrow}

\newcommand{\ta}{\tilde{a}}
\newcommand{\tb}{\tilde{b}}
\newcommand{\tc}{\tilde{c}}
\newcommand{\td}{\tilde{d}}
\newcommand{\te}{\tilde{e}}
\newcommand{\tf}{\tilde{f}}
\newcommand{\tg}{\tilde{g}}
\newcommand{\tih}{\tilde{h}}
\newcommand{\ti}{\tilde{i}}
\newcommand{\tj}{\tilde{j}}
\newcommand{\tk}{\tilde{k}}
\newcommand{\tl}{\tilde{l}}
\newcommand{\tm}{\tilde{m}}
\newcommand{\tn}{\tilde{n}}
\newcommand{\tio}{\tilde{o}}
\newcommand{\tp}{\tilde{p}}
\newcommand{\tq}{\tilde{q}}
\newcommand{\tr}{\tilde{r}}
\newcommand{\ts}{\tilde{s}}
\newcommand{\tit}{\tilde{t}}
\newcommand{\tu}{\tilde{u}}
\newcommand{\tv}{\tilde{v}}
\newcommand{\tw}{\tilde{w}}
\newcommand{\tx}{\tilde{x}}
\newcommand{\ty}{\tilde{y}}
\newcommand{\tz}{\tilde{z}}
\newcommand{\tA}{\tilde{A}}
\newcommand{\tB}{\tilde{B}}
\newcommand{\tC}{\tilde{C}}
\newcommand{\tD}{\tilde{D}}
\newcommand{\tE}{\tilde{E}}
\newcommand{\tF}{\tilde{F}}
\newcommand{\tG}{\tilde{G}}
\newcommand{\tH}{\tilde{H}}
\newcommand{\tI}{\tilde{I}}
\newcommand{\tJ}{\tilde{J}}
\newcommand{\tK}{\tilde{K}}
\newcommand{\tL}{\tilde{L}}
\newcommand{\tM}{\tilde{M}}
\newcommand{\tN}{\tilde{N}}
\newcommand{\tO}{\tilde{O}}
\newcommand{\tP}{\tilde{P}}
\newcommand{\tQ}{\tilde{Q}}
\newcommand{\tR}{\tilde{R}}
\newcommand{\tS}{\tilde{S}}
\newcommand{\tT}{\tilde{T}}
\newcommand{\tU}{\tilde{U}}
\newcommand{\tV}{\tilde{V}}
\newcommand{\tW}{\tilde{W}}
\newcommand{\tX}{\tilde{X}}
\newcommand{\tY}{\tilde{Y}}
\newcommand{\tZ}{\tilde{Z}}

\newcommand{\fl}{f_{\grl}}
\newcommand{\el}{E_{\grl}}
\newcommand{\sel}{\se_{\grl}}
\newcommand{\szl}{\sz_{\grl}}
\newcommand{\zl}{Z_{\grl}}
\newcommand{\famfl}{\{ \fl \} }
\newcommand{\sT}{{\cal T}}
\newcommand{\ab}{{\bf a}}
\newcommand{\Def}{{\bf Definition}}
\newcommand{\Pf}{{\bf Proof}}
\newcommand{\uG}{\bar{{\cal G}}}
\newcommand{\sG}{{\cal G}}
\newcommand{\sF}{{\cal F}}
\newcommand{\stG}{\tilde{{\cal G}}}
\newcommand{\tmu}{\tilde{\mu}}
\newcommand{\tpi}{\tilde{\pi}}
\newcommand{\tgrl}{\tilde{\grl}}
\newcommand{\sD}{{\cal D}}
\newcommand{\sDh}{\hat{\sD}}
\newcommand{\sP}{{\cal P}}
\newcommand{\sM}{{\cal M}}
\newcommand{\bsM}{\bar{{\cal M}}}
\newcommand{\sV}{{\cal V}}
\newcommand{\sS}{{\cal S}}
\newcommand{\sI}{{\cal I}}
\newcommand{\sW}{{\cal W}}
\newcommand{\sCZ}{{\cal C}{\cal Z}}
\newcommand{\tsT}{\tilde{\sT}}
\newcommand{\tsD}{\tilde{\sD}}
\newcommand{\tgre}{\tilde{\gre}}
\newcommand{\sH}{{\cal H}}
\newcommand{\tsH}{\tilde{\sH}}
\newcommand{\tsG}{\tilde{\sG}}
\newcommand{\tDel}{\tilde{\Grd}}
\newcommand{\sA}{{\cal A}}
\newcommand{\sB}{{\cal B}}
\newcommand{\sCF}{{\cal CF}}
\newcommand{\bsH}{\bar{\sH}}
\newcommand{\ntwob}{\bar{\nu_2}}
\newcommand{\nb}{\bar{\nu}}
\newcommand{\bF}{\bar{F}}
\newcommand{\bG}{\bar{G}}
\newcommand{\bI}{\bar{I}}
\newcommand{\bZ}{\bar{Z}}
\newcommand{\bh}{\bar{h}}
\newcommand{\bT}{\bar{T}}
\newcommand{\bx}{\bar{x}}
\newcommand{\bxi}{\bar{\xi}}
\newcommand{\bu}{\bar{u}}
\newcommand{\bg}{\bar{g}}
\newcommand{\bN}{\bar{\N}}
\newcommand{\bGrs}{\bar{\Grs}}
\newcommand{\bgrs}{\bar{\grs}}
\newcommand{\bmu}{\bar{\mu}}
\newcommand{\db}{\bar{d}}
\newcommand{\dx}{\dot{x}}
\newcommand{\n}{\norm}
\newcommand{\grfh}{\hat{\grf}}

\newcommand{\Lpfk}{\mbox{$\cap \hspace{-.44cm}\mid$}}
\newcommand{\Emn}{ R_{i_{-m } \ldots i_{-1} , i_0 \ldots i_{n-1}}}
\newcommand{\Rmn}{ R_{i_{-m } \ldots i_{-1} , i_0 \ldots i_{n-1}}}
\newcommand{\En}{ E_{ i_0 \ldots i_{n-1}}}
\newcommand{\Ep}{ E_{ i_0 \ldots i_{p-1}}}
\newcommand{\Sp}{ S_{ i_{0} \ldots i_{p-1}}}
\newcommand{\Sm}{ S_{ {i_{-m } \ldots  i_{-1} }}}
\newcommand{\Sn}{ S_{ i_{0} \ldots i_{n-1}}}
\newcommand{\Gum}{ \Gamma^u_{i_{-m } \ldots  i_{-1} }}
\newcommand{\Gsn}{ \Gamma^s_{ i_0 \ldots i_{n-1}}}
\newcommand{\Em}{ E_{i_{-m } \ldots i_{-1}}}
\newcommand{\Rpn}{ R_{i_{1 } \ldots i_{p} , n}}
\newcommand{\Rkn}{ R_{i_{1 } \ldots i_{k} , n}}
\newcommand{\Ran}{ R_{\bar{\alpha} , n}}
\newcommand{\Enewn}{ E_{ i_1 \ldots i_{n}}}
\newcommand{\Enewp}{ E_{ i_1 \ldots i_{p}}}
\newcommand{\Snewp}{ S_{ i_{1} \ldots i_{p}}}
\newcommand{\Snewn}{ S_{ i_{1} \ldots i_{n}}}
\newcommand{\Snewk}{ S_{ i_{1} \ldots i_{k}}}
\newcommand{\Snewjp}{ S_{ j_{1} \ldots j_{p}}}
\newcommand{\Rnewn}{ R_{i_{1 } \ldots i_{n-1} , i_{n}}}
\newcommand{\Sa}{ S_{\bar{\alpha} }}
\newcommand{\Can}{ C_{\bar{\alpha} , n}}
\newcommand{\Rbm}{ R_{\bar{\beta} , m}}
\newcommand{\Cbm}{ C_{\bar{\beta} , m}}
\newcommand{\Rak}{ R_{\bar{\alpha} , k}}
\newcommand{\fnewn}{ f_{ i_{1} \ldots i_{n}}}
\newcommand{\fnewp}{ f_{ i_{1} \ldots i_{p}}}

\def\bibsameauth{\leavevmode\vrule height .1ex depth 0pt width
2.3em\relax\,}

\begin{document}

\bibliographystyle{plain}

\title{Countable Markov partitions suitable for thermodynamic formalism}
\author{Michael Jakobson \thanks{University of Maryland,   Department of Mathematics, College Park, MD, USA} $\>$, Lucia  D. Simonelli \thanks{Abdus Salam International Centre for Theoretical Physics, Trieste, Italy}} 
\date{}


\maketitle

 {\centering To the memory of Roy Adler  \par } 
  
\begin{abstract}
\noindent  We study hyperbolic attractors of some dynamical systems with apriori given countable Markov partitions. Assuming that contraction is stronger than expansion we construct  new Markov rectangles such that their crossections by  unstable manifolds are Cantor sets  of positive Lebesgue measure. Using new Markov partitions we develop thermodynamical formalism and
  prove  exponential decay of
correlations and related properties for certain H\"older functions. The results are based on the methods developed by Sarig \cite{Sarig-1} -  \cite{Sarig-4}.

\end{abstract}
\section{Introduction} 

\noindent Examples of one-dimensional maps with countable Markov partitions go back to the Gauss transformation, and
further developments appeared in particular in \cite{Renyi-57}, \cite{Adler-73}, \cite{Walters-78}. Beginning in
\cite{Adler-79}, theorems about ergodic properties of such maps are often referred to as {\it Folklore}.\\

\noindent More recently an interest in such maps was motivated by works on
 ergodic and statistical properties of quadratic-like and H\'enon-like maps. The study of such maps is typically based on various tower constructions, see in particular \cite{J0},   \cite{J1}, \cite{BC}, \cite{BY},   \cite{WY}. 
Power maps defined on the tower satisfy hyperbolicity and distortion estimates.\\

\noindent Although the following remark is not directly related to the results of our paper, it might be useful for the further work 
in the area under discussion.\\
\begin{Rem} Since numerical evidence for the existence of a strange attractor was presented in the original H\'enon paper
\cite{Henon}, rigorous results were only obtained in unspecified small neighborhoods of one-dimensional maps.\\
\noindent One possible approach to that problem is to prove that  in a sufficiently small neighborhood of the classical H\'enon values, $a_H=1.4, b_H=.3$,
there is a positive Lebesgue measure set $M$, such that for $(a,b)\in M$, H\'enon maps $f_{a,b}$ have SRB measures with
strong mixing properties.\\
\noindent The main difficulty in that direction is to design a set of checkable numeric estimates which can be maintained through the
induction. In the one-dimensional case such estimates were used in \cite{LT}, \cite{YRH}, \cite{GLP}.\\
\end{Rem}
\noindent We study apriori given two-dimensional systems with countable Markov partitions satisfying hyperbolicity
and distortion conditions. In \cite{Jakobson} we proved strong mixing properties of such systems 
 assuming distortion condition D2, requiring boundedness of the quotient of the second derivatives over the first derivative. This condition is too strong for our purposes since it does not hold for power maps induced by quadratic and H\'enon maps. \\

\noindent Here we return to the more general setting of \cite{Jak-New-1}, \cite{Jak-New-2} where only
 boundedness of the quotient of the second derivatives over the square of first derivatives is assumed.
In order to study the decay of correlations  we require additionally that contraction in our models grows faster than expansion, see
condition H5 below.  Condition H5  naturally holds  for power maps generated by H\'enon-like maps with Jacobian less than $1$. We also require distortions of our  initial maps to be uniformly bounded, see condition BIV below. That is a standard requirement for maps defined on a tower.\\
\noindent The main new idea of the paper is to develop thermodynamic formalism by using special Markov rectangles  such that their intersections with unstable manifolds are Cantor sets of  positive Lebesgue measure.

\section{ Description of models and statement of results.}

\be

\item \textbf{Description of the model.} \\

\noindent The setting for our model is the same as in \cite{Jakobson} -  \cite{Jak-New-2}. To summarize, let $Q$ be the unit square.  Let $\xi = \{E_1, E_2, \ldots \} $ 
be a countable collection
of full-height, closed, curvilinear rectangles in $Q$. Hyperbolicity conditions that we will recall below imply 
 that the left and right boundaries of $E_i$ are graphs of
smooth functions $x^{(i)}(y)$ with $\left|{dx^{(i)}\over dy}\right| \leq
\gra$ for $0 < \gra < 1$. 

\noindent Assume that each $E_i$
lies inside a domain of definition of a $C^2$ diffeomorphism $F $
which maps $ E_i $ onto its image $S_i \subset Q$. The images $F|_{E_{i}}(E_{i})=f_i(E_i)=S_i$ are disjoint, full-width strips of $Q$ which are bounded from above
and below by the graphs of smooth functions $y^{i}(x),
\norm{{dy^{(i)}\over dx}} \leq \gra $.

We recall  geometric  and hyperbolicity conditions from \cite{Jak-New-2}.

 \item \textbf{Geometric conditions.} \\

For $z \in Q$, let $\ell_z$ be the horizontal line through $z$. We
define $\grd_z(E_i) = diam(\ell_z \bigcap E_i)$, $\grd_{i,max} = \max_{z
\in Q} \grd_z(E_i)$, $\grd_{i,min} = \min_{z \in Q} \grd_z(E_i)$. \\

\be
\item[G1.]  For $i \neq j$ holds  $int \ E_i \cap int \ E_j =\emptyset$, $int \ S_i \cap int \ S_j =\emptyset$ .
\item[G2.] $mes(Q \setminus \cup_i \  int \, E_i)  = 0$ where $mes$
stands for Lebesgue measure.
\item[G3.] for some $0<a \le b < 1$  and some ${\tC} \ge 1$ it holds that \[  {\tC}^{-1}a^i \le \grd_{i,min} \le    \grd_{i,max} \le {\tC}b^i. \] 
\ee

\item

{\bf Hyperbolicity conditions.} \\

Let $J_F(z)$ be the absolute value of the Jacobian determinant of $F$ at $z$.

There exist  constants $0 < \gra < 1$
and $K_0 > 1$ such that for each $i$   the map
$$ 
F(z)  = f_i(z) \mbox{ for } z \in   E_i
$$
satisfies 
\be
\item[H1.] $\n{F_{2x}(z)} + \gra \n{F_{2y}(z)} + \gra^2 \n{F_{1y}(z)}
\leq \gra \n{F_{1x}(z)} $
\item[H2.] $ \n{F_{1x}(z)} - \gra \n{F_{1y}(z)} \geq K_0$.
\item[H3.] $ \n{F_{1y}(z)} + \gra \n{F_{2y}(z)} + \gra^2 \n{F_{2x}(z)}
\leq \gra \n{F_{1x}(z)} $
\item[H4.] $ \n{F_{1x}(z)} - \gra \n{F_{2x}(z)} \geq J_F(z)K_0$.
\ee

We recall some notation.

Given a finite string $i_0 \ldots i_{n-1}$,  $n \ge 1$, we define inductively
\beq \label{fhE}
E_{i_0 \ldots i_{n-1}} = E_{i_0} \bigcap f_{i_0}^{-1}E_{i_1i_2 \ldots
i_{n-1}}
\eeq

Then, each set $E_{i_0 \ldots i_{n-1}}$ 
is a  full height subrectangle of $E_{i_0}$.  

Analogously, for a string
$i_{-m} \ldots i_{-1}$  we define

\[
{\Sm} = f_{ i_{-1}}(S_{i_{-m} \ldots i_{-2}} \bigcap
E_{ i_{-1}}) \] 
 
and get that  $S_{i_{-m} \ldots  i_{-1}}$ is a  full width strip in $Q$.
It is easy to see that $S_{i_{-m} \ldots  i_{-1}} = f_{ i_{-1}} \circ
f_{i_{-2}} \circ \ldots \circ f_{i_{-m}}(E_{i_{-m} \ldots  i_{-1}})$
and that $f_{i_0}^{-1}(S_{i_{-m} \ldots  i_{-1}i_{0}} )$ is a full-width substrip of
$E_{i_0}$. \\

We also define curvilinear rectangles ${\Rmn}$ by
\beq \label{oldR}
 {\Rmn} = {\Sm} \bigcap {\En} 
\eeq
If there are no negative indices then respective rectangle is  full height in $Q$.\\

The following  is a well known fact in hyperbolic theory, see \cite{Jak-New-2}.

\begin{Prop} \label{Prop1}
  Any $C^1$  map $F$   satisfying the above geometric
conditions G1--G3 and hyperbolicity conditions H1--H4 has a "topological attractor"
$$ \Lambda =\bigcup_{\ldots i_{-n} \ldots i_{-1}} \bigcap_{k\geq 1} 
S_{i_{-k} \ldots  i_{-1}}.$$ \vs

The infinite intersections
$$\bigcap_{k=1}^{\infty}S_{i_{-k} \ldots  i_{-1}}$$

define  $C^1$    curves
$y(x), \ | dy/dx| \leq  \gra,$  which are the unstable manifolds
for the points of the attractor.\\

The infinite intersections

$$\bigcap_{k=1}^{\infty} E_{i_0 \ldots i_{k-1}}$$

define  $C^1$
curves 
$x(y), \ | dx/dy| \leq  \gra,$  which are the stable manifolds
for the points of the attractor. \\

The infinite intersections

 \[ \bigcap_{m=1}^{\infty} \bigcap_{n=1}^{\infty} {\Rmn} \]
define
points of the attractor. 
\end{Prop}

\item {\bf Distortion condition. } \\ 
 
We  formulate certain assumptions on the second
derivatives.   We  use the 
distance function $d((x,y), (x_1,y_1)) = \max(\n{x - x_1}, \n{y -
y_1})$ associated with
the norm $\n{v} =
\max(\n{v_1}, \n{v_2})$ on vectors $v = (v_1, v_2)$. \\

Our first condition is the same as in \cite{Jak-New-2}; we recall it below.

Let $f_i(x,y) = (f_{i1}(x,y), f_{i2}(x,y))$. 
We use
$f_{ijx}, f_{ijy}, f_{ijxx}, f_{ijxy}$, etc. to denote partial derivatives 
of $f_{ij}, j = 1,2$. 

We define 

\[ \n{D^2f_i(z)} = \max_{j=1,2, (k,l) = (x,x), (x,y), (y,y)}
\n{f_{_{ijkl}}(z)}. \]

We assume that there exists a constant $C_0 > 0$ such that the following {\it Distortion Condition} holds:

\be
\item[D1.] $ \ \  \sup_{z \in E_i,\> i \geq 1}  \frac{\n{D^2f_i(z)}}
 { \norm{f_{i1x}(z)}}\grd_z(E_i) < C_0 $. 
\ee \vs

\item \textbf{A result about systems satisfying geometric, hyperbolicity and distortion conditions:} \\

Our conditions imply the following Theorem proved in \cite{Jak-New-1} and  \cite{Jak-New-2}.

\begin{theorem} \label{Theorem1}
 Let $F$ be a piecewise smooth mapping as above
satisfying the geometric conditions G1--G3, the hyperbolicity conditions
H1--H4, and  the distortion condition D1. 

Then, $F$ has an $SRB$ measure $\mu$ supported on 
$\Lambda $ whose basin has full Lebesgue measure in
$Q$.   The dynamical system   $(F,\mu)$ satisfies the following properties.\\
\be
\item $(F,\mu)$ is measure-theoretically isomorphic to a 
Bernoulli shift.
\item $F$ has finite entropy with respect to the measure
$\mu$,  and the entropy formula holds
\beq \label{ent_form1}
 h_{\mu}(F) = \int log|D^uF| d\mu 
\eeq
where $D^uF(z)$ is the norm of the derivative of $F$ in the
unstable direction at $z$. 
\item
\beq \label{ent_form}
 h_{\mu}(F) = \lim_{n \rarrow \infty} \frac{1}{n} \log \n{DF^n(z)} 
\eeq
where the latter limit exists for Lebesgue almost all $z$ and is
independent of such $z$. 

\ee
\end{theorem}

\item \textbf{Additional distortion and hyperbolicity conditions and statement of the main theorems.}\\
\be
\item Properties of the function  $\phi(z) = - \log(D^u F(z)) $ are important when applying thermodynamic formalism to hyperbolic attractors.
We consider systems
satisfying conditions of Theorem \ref{Theorem1}
and some extra
   hyperbolicity conditions, which   can be used for  power maps arising from Henon-type diffeomorphisms.\\

 We explore a general principle that can be stated as :
{\it Contraction  increases faster than expansion} -
see hyperbolicity condition H5  below.
 For such systems we construct new  Markov partitions such
that the pullback of  $\phi(z)$
 into a   respective symbolic space  is a  locally H{\" o}lder function. \\ 
 
 New  Markov rectangles are
Cantor sets, such that their one-dimensional crossections by $W^u(z)$ have positive Lebesgue measure. \\

\item We consider a class of system which satisfy conditions of the Theorem \ref{Theorem1} as well as the following additional assumptions.
\be
\item

{\bf Bounded Initial  Variation}. \vs 
\be
\item [BIV.] 
 There exists $B_0 >0$ such that for all $i$ and all 
 
 $$\{ z_1= (x_1,y_1),  z_2= (x_2,y_2)
 \} \in E_i
 $$ 
 
 holds 
\beq \label{BV}
\n{\log f_{i1x}(z_1) - \log f_{i1x}(z_2)} <B_0.
\eeq
\ee \vs
BIV does not allow unbounded oscillations of widths for initial rectangles;
it implies that $a=b$ in condition G3:
\beq \label{modifG3}
  {\tC}^{-1}a^i \le \grd_{i,min} \le    \grd_{i,max} \le {\tC} a^i 
\eeq
\item {\bf Contraction grows faster than expansion}. \\
We assume that there is a constant $a_1$ satisfying 
\beq \label{a1}
 0 < a_1 < a
\eeq
where $a$  is from (\ref{modifG3}),
such that for each $j$, for each $z \in E_j$, and 
for any vector $v$ in the stable cone $K^s_{\gra}(z)$ holds

\be \label {strongcontraction}
\item[H5.] $  \n{Df_j^{-1}v} \ge  a_1^{-j} \n{v}. $
\ee \vs
Condition H5 means that up to a uniform factor, contractions of $f_j$ grow faster than expansions.
In particular it implies that up to a uniform factor,
heights of
$ S_i \cap E_k$ are smaller than widths of $E_k$ for all $k \le i$. 
\ee
\begin{Rem} Uniform hyperbolicity and distortion conditions were used in \cite {Jak-New-1},  \cite {Jak-New-2} 
to extend the classical approach of \cite{Anosov-Sinai} and to study ergodic properties of systems with countable Markov
partitions. By combining D1 and H5 we can  add methods of thermodynamic formalism.
\end{Rem}
 
\item
Let ${\cal H}_{\gamma}$ be the space of functions on $Q$ satisfying 
H\"older property with exponent $\gamma$ 
\[ \norm{\phi(x)-\phi(y)} \le c\norm{x-y}^{\gamma}.\]

 We state our main Theorems.

\begin{theorem} \label{main} ({\bf Exponential Decay of Correlations})\\
 Let $F$ be a piecewise smooth mapping as above,
satisfying  geometric conditions G1-- G3,  hyperbolicity conditions
H1--H5,   distortion condition D1, and the BIV condition. Then the system $(F, \mu)$ has exponential decay of correlations  for
 $f \in \mathcal{H}_{\gamma}$ and $g \in L^{\infty}(\mu)$.
Namely there exists $0<\eta<1,  \  \eta=\eta(\gamma)$, such that
\beq  
\left | \int f \> (g \circ F^n) \> d\mu - \int f \> d\mu \int g \> d\mu \right | < C(f,g)\eta^{n}.
\eeq
\end{theorem} 

\begin{theorem} \label{CLT} ({\bf Central Limit Theorem})\\
Let $(F, \mu)$ satisfy the assumptions of Theorem \ref{main} and suppose that  $f \in \mathcal{H}_{\gamma}$. If $\int f d\mu = 0$ and $f$ cannot be expressed as $g-g \circ F$ for $g$ continuous, then there is a positive constant $\sigma = \sigma(f)$, such that for every $t \in \R$,

$$
\lim_{n \to \infty} \> \> \mu \left \{  x: \frac{1}{\sqrt{n}} \sum_{k=0}^{n-1}  f \circ F^{k}(x) \> < \> t  \right \}= \frac{1}{\sqrt{2 \pi \sigma^{2}}} \int_{-t}^{\infty}e^{\frac{-s^{2}}{2\sigma^{2}}} \> ds.
$$

\end{theorem}

\ee
\ee
\section{H{\" o}lder  properties of  $\log(D^u F(z)) $ and Markov partitions in the phase space. }
\be
\item
The key step toward the proof of Theorem \ref{main} is to establish that the pullback
of the function $\log  D^uF$ into the respective symbolic space  is H\"older continuous.
Then one can follow the Ruelle-Bowen approach (\cite{Ruelle-1}, \cite{Bowen}), in particular using results of Sarig \cite{Sarig-1} - \cite{Sarig-4} to develop thermodynamic formalism for the systems under consideration. H\"older  properties of the pullback
of  $\log  D^uF$ into  symbolic space follow from H\"older  properties 
of  $\log  D^uF$ in the phase space. \\
In order to get an appropriate symbolic space, we construct a partition of a subset  of positive measure $\cal{C} \subset$ $\Lambda$,
such that the first return map to $\cal{C}$  is Markov. Elements of the Markov partition of  $\Lambda$
are elements of the Markov partition of $\cal{C}$ and their orbits
before the first return. \\ 
Elements $C_i$ of the Markov partition $\cal{C}$  have the following property. For $z \in C_i$
the crossection $C_i^u(z)$ of $C_i$ by $W^u(z)$ is a Cantor set of
positive linear Lebesgue measure.
\item
Cantor sets that we construct are inscribed in curvilinear rectangles\\  ${\Emn} $. Recall that
 ${\Emn} $ are bounded from above and below by 
 arcs of  unstable curves $\Gamma^u$, which are
images of some pieces of the top and bottom of $Q$,
and  from the left and right by  arcs
 of   stable curves ${\Gamma^s}$ which are  preimages of some pieces of the left and right boundaries of $ Q$ . \\

For $x \in {\Rmn}\cap \Lambda$ let 
$$\Gamma^s(x,{\Rmn}) = W^s(x) \cap {\Emn}$$
$$\Gamma^u(x,{\Rmn}) =  W^u(x) \cap {\Emn}.$$
We define the height of ${\Emn}$ as
 \[ H({\Emn}) = sup_{x} \left |\Gamma^s(x,{\Rmn}) \right |\]
 The width $ W({\Emn})$ is defined similarly.\\

Hyperbolicity conditions imply that stable boundaries of rectangles belong to stable cones. Since standard horizontal lines belong to unstable cones, and stable and unstable cones are separated, we get the following. \\

For every $\gre_1$ there is an $\gre_0$ such that if

\beq \label {heightwidth} 
 H({\Emn}) < \gre_0 w_{min}({\Emn}) 
\eeq

then for all $m \ge 1$, $n\ge 1$, the ratio of lengths of any two unstable curves $\Gamma^u (x,{\Rmn})$
is bounded by $1 \pm \gre_1$. Similarly
it follows from hyperbolicity conditions, that if $l$ is the length of a
standard horizontal crossection of ${\En}$ through a point  $x \in {\Emn}$, then
for some $c_0 $,
\beq \label {crossection} 
  \norm { \Gamma^u(x,{\Rmn})  }< c_0 l. 
\eeq

\item {\bf Admissible objects}.\\

 We define the following  strings of indices as { \it admissible}:

\item [A1.] A string
$ \bar{\gra} = [i_{1} \ldots i_k ] $ is admissible if
for each $l= 1,2, \ldots k-1$ it holds that
\beq \label{A2new}
   \sum_{m=1}^{l}  i_m \ge i_{ l+1} 
\eeq

A rectangle ${\Rmn}= $ is admissible if  the string $ [i_{-m} \ldots i_{-1} i_0 \ldots i_{n-1} ] $ is admissible.
It follows from definition that if ${\Rmn} $, $m \ge 0, n\ge 1$  is admissible, then all rectangles obtained by moving the comma in the index to the left or to the right are admissible.
A one-sided sequence $i_1i_2 \ldots i_n \ldots$ is admissible if all strings $[i_1 \ldots i_n ]  $
are admissible. 
 
 In particular all strings $[ij], i\ge j$ are admissible, and thus, the respective rectangles are admissible.\\

Note that distortion estimates may be satisfied on non-admissible rectangles if their heights are small enough, but we ignore that possibility, and organize our construction based on condition A1. \\

\item
We estimate the variation of $\log  D^uF$
on  admissible  two-dimensional  curvilinear rectangles 
${\Emn}$. For any function $a(x,y)$   the variation of $a(x,y)$ over a rectangle $R$ is defined as
\beq \label {vardef}
 var(a(x,y)) | R = \sup_{(x_1,y_1) \in R,(x_2,y_2) \in R} \n{a(x_1,y_1)-a(x_2,y_2)}
\eeq

The function $\log  D^uF$ is {\bf locally  H{\" o}lder} on admissible rectangles
 ${\Emn}$ if for $m \ge 0$, $n \ge 1$
  the variation of $\log  D^uF$ on  ${\Emn}$ satisfies
\beq \label {var1}
 var(\log D^uF) | {\Emn}< C\theta_0^{min(m,n)}
\eeq
for some $C>0$, $\theta_0 < 1$.\\

Note that on initial rectangles $E_i$ estimate (\ref{var1}) is satisfied because of BIV.

\item The proof of the following Proposition is similar to the proof of Proposition 5.1 in \cite{Jakobson}.

\begin{Prop} \label {unifdist1}
For any admissible string $[i_{-m} \ldots i_{-1}i_{0} \ldots i_{n-1}]$ 
the variation of  $\log  D^uF$ on ${\Emn}$ satisfies (\ref{var1}) with some $C$ and $\theta_0$
independent of $m,n$ and determined by hyperbolicity and distortion conditions.
\end{Prop}

Note that   in \cite{Jakobson} and \cite{Jakobson1} we proved H\"older property of $\log  D^uF$ on arbitrary rectangles
${\Emn}$. Here we prove it for admissible rectangles. 
\vspace{0.2cm}

{\bf Proof}.

\be
\item Admissible rectangles ${\Emn}$ are bounded from above and below by some arcs of two unstable curves $\Gamma^{u}_{i_{-m}...i_{-1}}$ which are images of some pieces of the top and bottom of $Q$ and from left and right by some arcs of two stable curves $\Gamma^{s}_{i_{0}...i_{n-1}}$ which are preimages of some pieces of the left and right boundaries of $Q$.

\vspace{0.2cm}

\noindent Let $z_{1}, z_{2} \in \Emn \cap \Lambda$. We consider two points $z_{3}, z_{4}$ such that $W^{s}(z_{3})=W^{s}(z_{4})$ and for which we can connect $z_{1}$ to $z_{3}$ and $z_{2}$ to $z_{4}$ along their respective unstable manifolds. We define the following curves inside $\Emn$,
$$
\gamma_{1}=\gamma(z_{1},z_{3}) \subset W^{u}(z_{1})
$$
$$
\gamma_{2}=\gamma(z_{2},z_{4}) \subset W^{u}(z_{2})
$$
$$
\gamma_{3}=\gamma(z_{3},z_{4}) \subset W^{s}(z_{3}).
$$
\noindent Now we bound each term on the right hand side of the inequality
\bqn \label{p11}
|logD^{u}F(z_{1}) - logD^{u}F(z_{2})| & \leq & \nonumber \\
|logD^{u}F(z_{1}) - logD^{u}F(z_{3})|  &+&  \nonumber \\
 |logD^{u}F(z_{3}) - logD^{u}F(z_{4})| &+&  \nonumber \\
 |logD^{u}F(z_{4}) - logD^{u}F(z_{2})| \nonumber \\
\eqn
Estimates of $|logD^{u}F(z_{1}) - logD^{u}F(z_{3})|$ and 
 $|logD^{u}F(z_{4}) - logD^{u}F(z_{2})|$ are the same as estimates $(15)-(28)$ in the proof of Proposition 5.1 in 
\cite{Jakobson}. Then we get

\beq \label{p15}
|logD^{u}F(z_{1}) - logD^{u}F(z_{3})| < C_{2}\frac{1}{K_{0}^{n}}.
\eeq

\noindent Similarly,

\beq \label{p16}
|logD^{u}F(z_{2}) - logD^{u}F(z_{4})| < C_{2}\frac{1}{K_{0}^{n}}.
\eeq

\item The second part of the proof, depending on $m$, follows again the ideas in \cite{Jakobson} and \cite{Jakobson1} but also utilizes condition H5. We are left with estimating the difference
\beq \label{p21}
|logD^{u}F(z_{3}) - logD^{u}F(z_{4})|.
\eeq
\noindent Note that the BIV condition implies that the above difference is uniformly bounded on full-height rectangles. From \cite{Jak-New-2} we get that the hyperbolicity conditions imply that any unit vector in $K^{u}_{\alpha}$ at a point $z \in E_{i}$, in particular a tangent vector to $W^{u}(z)$, has coordinates $(1,a_{z})$ with $|a_{z}|<\alpha$. Thus we need to estimate

\beq \label{p22}
log|F_{1x}(z_{3})+a_{z_{3}}F_{1y}(z_{3})|-log|{F_{1x}(z_{4})-a_{z_{4}}F_{1y}(z_{4})}|.
\eeq
\noindent Now we are moving along $\gamma_{3} \subset W^{s}(z_{3})$ connecting $z_{3}$ and $z_{4}$. \\ We cover $\gamma_{3}$ by rectangles $\tilde{R_{k}}$ for which the widths $\Delta x$ and lengths $\Delta y$ satisfy $|\Delta x| < \alpha |\Delta y|$. As in \cite{Jakobson}, we  get (\ref{p22})  
by estimating differences

\beq \label{p23}
|logF_{1x}(\tilde{z})-logF_{1x}(\tilde{z}')|
\eeq

\noindent for $\tilde{z}, \tilde{z}' \in \tilde{R}_{k} \cap W^{s}(z_{3})$, and

\beq \label{p24}
|a_{z_{3}}-a_{z_{4}}|.
\eeq

\noindent 
 To estimate (\ref{p23}) we use the mean value theorem and get on each rectangle $\tilde{R}_{k}$ an estimate not exceeding

\beq \label{p25}
Const \cdot \sup_{z\in \tilde{R}_{k}} \frac{|f_{1ij}(z)|}{|f_{1x}(z)|}\Delta y.
\eeq
 Let $\Gamma_{3}=W^{s}(z_{3}) \cap R_{i_{-1},i_{0}}$.
 The sum of the contributions  from (\ref{p25}) is bounded by

\beq \label{p251}
Const \sup_{z\in R_{i_{-1},i_{0}}} \frac{|f_{1ij}(z)|}{|f_{1x}(z)|} |\gamma_{3}|.
\eeq

 Since the distortion condition D1 is expressed using the width of $E_{i_0}$, we  use 
 condition H5. \\
 
 On an admissible rectangle $\Emn$, the string  $[i_{-m} \ldots i_{-1}i_{0}]$ is
 admissible. Let $h_{max}$ be the 
maximal height of $S_{i_{-m} \ldots i_{-1}}$
and let $w_{min}$ be the minimal width of $E_{i_0}$. From condition H5 we know that contraction of $f_i$ is stronger than $a_1^i$.
Since the rectangle is admissible, the sum of the indices satisfies $i_{-m}+ \ldots i_{-1} \ge i_0$. Since contraction of the composition is stronger than $a_1^{i_0}$, it follows that
\beq \label {hw2new}
\frac {h_{max}}{  w_{min} } < C \Big ( \frac{a_1}{a} \Big )^{i_0}.
\eeq

 As each index is at least $1$ we get that $i_0 \ge m$ and thus

\beq \label {hw2}
 h_{max} < \Big ( \frac{a_1}{a} \Big )^{m}{\tC}^{-1} w_{min}. 
\eeq

Therefore
 the heights of rectangles $R_{i_{-m} \ldots i_{-1},i_0}$ decay exponentially
comparatively to the width
of $E_{i_0}$. Since $|\gamma_{3}|<h_{max}$ and $w_{min}<\delta_{z}(E_{i_{0}})$ for any $z \in E_{i_{0}}$,  we can apply D1 and  obtain the following  bound for 
the sum of the contributions from (\ref{p23}),
\beq \label{p253}
C_{3} \Big ( \frac{a_1}{a} \Big )^{m}
\eeq
\vspace{0.2cm}
\noindent We estimate  (\ref{p24}) as in \cite{Jakobson}, \cite{Jakobson1}. We assume by induction

\beq \label{p24a}
|a_{z_{3}}-a_{z_{4}}|< c_1\theta_1^m
\eeq

As in \cite{Jakobson1} one can assume by taking if needed instead of $F$ some power of $F$ 
\beq \label{K0alpha}
\frac{1}{K_0^2} + \alpha^2 < 1
\eeq
Note that differently from  \cite{Jakobson1} the variation of $\log D^u F$ along stable
manifolds inside admissible rectangles is controlled not by using bounded distortions
of inverse maps, but from (\ref {hw2}) and (\ref{K0alpha}).\\ 

With that modification we prove like in Lemma 5.2 from \cite{Jakobson}

\beq \label{p27}
|a_{F(z_{3})}-a_{F(z_{4})}| < c_{1} \theta_{1}^{m+1}.
\eeq
where 
\beq \label{K0alphaa1a}
max \{ \frac{1}{K_0^2} + \alpha^2, \frac{a_1}{a} \} < \theta_1 < 1
\eeq
Combining (\ref{p253}) and (\ref{p27}) gives us 
\beq \label{p28}
|logD^{u}F(z_{3})-logD^{u}F(z_{4})|<C_{4}\theta_{1}^{m}.
\eeq
\noindent Finally 
combining (\ref{p15}), (\ref{p16}), and (\ref{p28}) concludes the proof of  Proposition \ref{unifdist1},
if we take $\theta_0 < 1$ satisfying

\beq \label{theta0theta1}
\theta_0 > max \{ \frac{1}{K_0},  \theta_1 \}.
\eeq

\ee

%

\item {\bf Construction of full height Cantor sets}. \\

We define full height Cantor sets $C_n$ inside full height rectangles $E_n$, $ n \ge 1$ by
\beq \label {Cantor2}
C_n= \bigcap_{m=1}^{\infty} \bigcup_{k =1}^m E_{i_{0}  \ldots i_{k-1}}
\eeq
where $i_0=n$, $k \ge 1$ and all $[i_{0} \ldots i_{k-1}]$ are admissible strings. \\

As an example let us consider several rectangles with indices starting from $1$.
It follows from the definition that $[11]$ is admissible and  $[1i], i>1$ are not. Then $R_ {[11]}$ is the only
 defining rectangle of order two for the Cantor set  $C_1$ and that the other $R_ {[1i]}$ are gaps.  The next defining rectangles of order three for $C_1$ are
$R_ {[111]}$ and $R_ {[112]}$, of order four $R_ {[1111]}$,  $R_ {[1112]}$,  $R_ {[1113]}$,
 $R_ {[1121]}$,  $R_ {[1122]}$,  $R_ {[1123]}$,  $R_ {[1124]}$ and so on. \\

 As each index is at least $1$ we get from the definition of admissible rectangles  that inadmissible indices 
satisfy $i_N > N$. Geometric condition G3 implies that on any unstable manifold, the relative measure 
of the union of rectangles with indices greater than $N$ decays exponentially.\\
Then uniformly bounded distortion implies that 
the total relative linear measure of gaps in an unstable manifold of any defining  rectangle of order $N$
has uniform exponential decay. \\

This implies the following Corollary.

\begin{Cor} \label{positivemeasure}
{\it There is a $c_0>0$ such that for any initial full height rectangle $E_n$ and respective 
 full height Cantor set $C_n$ constructed inside $E_n$ and for any
$z \in C_n \cap \Lambda$ 
the relative linear measure of $C_n$ in $ W^u(z,E_n)$ is greater than $c_0$.
Moreover that relative measure tends to one when $n \rar \infty$.}
\end{Cor}
 Also as in \cite{Jakobson},  Remark 5.10,  we get the following Corollary from distortion estimates.
\begin{Cor} \label{constwidth}
{\it Let ${\En}$ be a full height admissible rectangle of order $n$. \\
Then for any two points  $z_1, \> z_2 \in {\En} \cap \Lambda$, it holds that}
{\vspace{.1in}}
\beq  \label {bndratiowidths}
\frac{\norm{ W^u(z_1, {\En})} } { \norm{ W^u(z_2, {\En})}} < c
\eeq
\end{Cor}
As Corollary \ref{constwidth} is valid for all defining rectangles we get
\begin{Cor} \label{bndratiolinmeasures}
{\it Let $C_n$ be the full height Cantor set constructed inside  $E_n$,
let $z \in C_n \cap \Lambda$, and let $\norm{ W^u(z, C_n)}$ be the linear Lebesgue measure of $C_n$ in $W^u(z, E_n)$. Then for any two points  $z_1, \ \ z_2 \in C_n \cap \Lambda$, it holds that }
{\vspace{.1in}}
\beq  \label {bndratiowidths1n}
\frac{\norm{ W^u(z_1, C_n)} } { \norm{ W^u(z_2, C_n)}} < c
\eeq
\end{Cor}
\item { \bf Markov properties of $C_n$}. \\

Every Cantor set $C_i$ is determined by its defining rectangles and equivalently by its gaps. Defining rectangles are labeled by admissible strings $[i_1i_2 \ldots i_n]$ satisfying A1.
Gaps are labeled by nonadmissible strings $[i_1i_2 \ldots i_nj]$, where $[i_1i_2 \ldots i_n]$
is admissible and $j > i_1+i_2 + \ldots + i_n$.\\

For example gaps of $C_3$ are: \\

 $E_{3k}, \>  k >3$ - gaps of order 1,\\
 
$E_{31k}, \> k>4$, $E_{32k}, \> k>5$, $E_{33k}, \> k>6$ - gaps of order 2,\\

$E_{311k}, \> k>5$ - gaps of order 3, and so on. \\

\noindent The following example illustrates Markov relations between Cantor sets $C_i$. \\

\noindent Consider an admissible rectangle $E_{112}$. Then $F(E_{112}) = R_{1,12}$ is a subrectangle
of a non-admissible rectangle $E_{12}$. However, because of H5, $F^2(E_{112}) = R_{11,2}$ has height less than the width of $E_2$. \\

We will use the notation

\[
C_{112} = C_1 \cap E_{112}
\]

and in general,

\[
C_{i_1 \ldots i_s}= C_{i_1} \cap E_{i_1 \ldots i_s},
\]
\[
C_{i_1 \ldots i_k, i_{k+1} \ldots i_l}= C_{i_{k+1} \ldots i_l} \cap R_ {i_1 \ldots i_k, i_{k+1} \ldots i_l}.
\]

We have that
\beq  \label {markovcover}
F^2(C_{112})  \supset C_{11,2}.
\eeq
As the sum of indices in $[112]$ is greater than $2$, the inclusion in (\ref{markovcover}) is not
an equality. Namely strings $[23]$ and $ [24]$ are not admissible, but $[1123]$ and $[1124]$ are
admissible. The image $F^2(C_{112})$ covers respective slice of $C_2$ and also covers some parts of 
gaps $C_{23}$ and $C_{24}$ used in the construction of the Cantor set $C_2$. \\

Consider the union of all  full height Cantor sets $C_n$ :
\beq  \label {Cantorunion}
{\cal C} = \bigcup_{n=1}^{\infty} C_n.
\eeq
For $z \in C_{i_1 \ldots i_n}$ let us denote $W^u(z,C_{i_1 \ldots i_n})= C_{i_1 \ldots i_n} \cap
W^u(z,E_{i_1 \ldots i_n})$. \\

Let $T$ be the first return map to $\cal C$ generated by $F$ . Similarly to (\ref{markovcover}) above we get the following Markov properties.
\begin{Prop} \label{Markovprop}
If $[i_1 \ldots i_n]$ is admissible, $z \in C_{i_1 \ldots i_n}$, $T(z) \in C_{i_n}$, then
\beq  \label {MP}
T(W^u(z,C_{i_1 \ldots i_n})) \supset W^u(T(z), C_{i_n}).
\eeq
\beq  \label {MP1}
 T(W^s(z,C_{i_1 \ldots i_n})) \subset W^s(T(z),C_{i_1 \ldots i_{n-1}, i_n})
\eeq

\end{Prop}
\item {\bf Estimates of measure}. \\

Let $T_1$ be the first return map to $C_1 \subset {\cal C}$. Next we  estimate the measure of points in $C_1$
 which return to $C_1$ after at least  $n$ iterates of $F$.

\begin{Prop} \label{FR}
Let $B_n$ be the  set of points in the domain of $T_1$ such that the return time for
$z \in B_n$ is greater than $n$. For some $C >0$ and $0< \beta <1$,
\beq  \label {estim1}
\mu(B_n) < C\beta^n.
\eeq
\end{Prop}

We begin by proving (\ref{estim1}) for the first return map $T$ onto $\cal C$. \\

\noindent Suppose $x=(x_0x_1 \ldots  x_n \ldots ) \in {\cal C}$, $y= Fx =(x_1x_2 \ldots x_n x_{n+1} \ldots)
\notin {\cal C}$, and   $F^nx  \in {\cal C}$ is the first return.\\

\noindent As $y \notin {\cal C}$, there is a coordinate $y_k$ such that 
$y_k = x_{k+1} > y_0 + \ldots +y_{k-1} = x_1 + \ldots + x_k $. If  $k \ge n$,
then
$$  x_{k+1} > x_1 + \ldots x_k \ge  x_n + \ldots + x_k $$
which contradicts that $F^n(x) \in {\cal C}$. \\

\noindent Suppose $k < n$ satisfies 
$$  x_1 + \ldots + x_k <  x_{k+1} $$
Then $x_{k+1} \ge k+1 $ because each coordinates is at least $1$. \\

As all images $F^k x , \  1 \le k < n$ do not belong to $\cal C$,
 there is a coordinate $x_N, \ N\ge k$ such that 
$$  x_k + \ldots + x_N <  x_{N+1}. $$
If $N\ge n$, then
$$  x_n + \ldots + x_N <  x_{N+1} $$
and we get a contradiction as above to $F^nx  \in {\cal C}$. So $N < n$ and we get
\beq \label{estimxn}
N+1 = 1+(k+1)+(N-k-1) \le x_k +x_{k+1}  + \ldots +x_{N}  < x_{N+1}
\eeq
 Proceeding as above we get
\beq  \label {atleastn}
x_n \ge n.
\eeq
As widths of $E_i$ decay exponentially, the measure of the collection of $x$ satisfying (\ref{atleastn}) decays exponentially. Since the measure of ${\cal C}$ is positive, we get that the measure of points which do not return to $\cal C$
after $n$ iterates is less than
\beq  \label {estim11}
 C_1\beta_1^n.
\eeq
 for some $C_1 >0$ and $0< \beta_1 <1$.\\
 
Next we note that
if $z \in C_1$ is mapped by $T$ into $C_i$, then because all transitions from $C_i$ to $C_1$ are admissible,
$z$ will be mapped into $C_1$ by the next iterate of $F$. Therefore points which do not return into
$C_1$ after $n$ iterates are subdivided into two subsets: points which did not return into $\cal C$
after $n$ iterates and points which returned into $C_i , \  i> 1$, after $n-1$ iterates and at the next iterate
were not mapped into $C_1$. Because of uniformly bounded distortion the measure of the second set
is less than
\beq  \label {estim111}
 C_1\beta_1^{n-1}(1-\gamma_0)
\eeq
where $0< \gamma_0 < 1$. Then  (\ref{estim1}) follows from (\ref{estim11}) and (\ref{estim111}) if we take
$C = 2C_1$ and $\beta= max \{ \beta_1, 1- \gamma_0 \}$.\\

This concludes the proof of Proposition \ref{FR}.

\item {\bf First return maps and respective transition matrices. } \\

Note that on every unstable leaf, relative measures of $C_i$ inside $E_i$ are uniformly bounded away from $0$. Together with uniformly bounded distortion it implies that in the orbits of the first return map
each gap is substituted by a union of new gaps and a Cantor set of relative measure greater than some
uniform $c_0 >0$. Then at the end of that construction we get, up to a set of measure zero,
\beq  \label {CantorunionCantor}
C_i  = \bigcup_{k} \{ T^{-1}(TC_i \cap C_k) = C_{ik} \}
\eeq
Note that $C_{ik}$ can be unions of several Cantor sets which belong to disjoint full height rectangles.
For example  consider $C_{1113} \subset E_{1113}$. Admissible rectangle $E_{1113}$  is mapped as follows,
$$E_{1113} \to E_{113} \to E_{13} \to E_{3} $$
As $E_{113}$, $E_{13}$ are inadmissible, we get that $T= F^3 $ maps $C_{1113}$  onto
$C_3$ in a Markov way.
Similarly $T= F^3 $ maps $C_{1123}$  onto
$C_3$ in a Markov way. \\
The correct labeling is provided by respective strings of the original alphabet starting width $i$ and ending with $k$.\\

To get an authentic Markov partition which generates a transition matrix of $0$'s and $1$'s we partition each
$C_{j_0j_1}$ into subsets
\beq  \label {strings1}
C_{i_0i_1 \ldots i_{n-1}}
\eeq
where 
\beq  \label {strings2}
i_0=j_0, i_1,\ldots ,i_{n-1}=j_1
\eeq
is admissible, and for all $k>0$
\beq  \label {strings3}
 i_k,\ldots, i_{n-1}
\eeq
are not admissible.\\
In other words the first return map $T$ maps  $C_{i_0 i_1 \ldots i_{n-1}}$ onto a full width substrip of $C_{ i_{n-1}}$
and all intermediate images of $C_{i_0 i_1 \ldots i_{n-1}}$ belong to various gaps.\\

As in the proof of (\ref {estim1})  we get that the length of the above strings starting from $j_0$ and ending
with $j_1$ is at most $j_1-j_0+2$. \\
The union of  $C_{j_0 i_1 \ldots i_{n-2}j_1}$ forms a Markov partition

\beq  \label {MP11} 
 { \cal MP} = \{C_{j_0i_1 \ldots i_{n-2}j_1}\}
\eeq
   of $\cal C$ . Using Proposition \ref{Markovprop} we get
\begin{Lem} \label{MarkovPartition1}
To any one-sided $T$-admissible sequence of transitions
 $$C_{j_0i_1^1 \ldots i_{n_1-2}^1j_1} \rar C_{j_1   i_1^2\ldots i_{n_2-1}^2  j_2} \rar \ldots $$
corresponds a unique one-sided sequence of the original alphabet
\beq \label{os}
j_0, \> i_1^1,\ldots, i_{n_1-1}^1=j_1, \> i_1^2,\ldots, i_{n_2-1}^2= \ldots
\eeq
such that 
$$C_{j_0i_1^{1} \ldots i_{n_1-2}^{1}j_1} \cap T^{-1}( C_{j_1 i_1^2\ldots i_{n_2-1}^2  j_2}   )  \cap \ldots $$
coincides with the stable manifold labeled by (\ref{os}).\\

The union of elements of  ${\cal MP}$ and all intermediate iterates of $C_{j_0i_1 \ldots i_{n-2}j_1}$ form a tower over
 $\cal C$. 
Elements of this tower form a Markov partition of  the attractor $\Lambda$.\\

Up to a set of $\mu$ measure zero, any point of the attractor is uniquely labeled  by a two-sided sequence
of admissible transitions 

\beq \label {admtwo}
 \ldots \rar  C_{j_{-1}i_1^{-1} \ldots i_{n_{-1}-2}^{-1}      j_0} \rar C_{j_0i_1^{1} \ldots i_{n_1-2}^{1}j_1}  \rar C_{j_1 i_1^2,\ldots i_{n_2-2}^2  j_2} \rar \ldots 
\eeq
\end{Lem}
\item 

We consider a new
alphabet $\Omega$ corresponding to the elements of the tower from Lemma \ref{MarkovPartition1} and get a subshift 
\beq \label {newsubshift}
(\Omega,X,\sigma)
\eeq

Recall that a subshift is  topologically mixing if for any states $a$ and $b$ there is $n(a,b)$ such that for $n \ge n(a,b)$ there is an admissible word of length $n$
starting from $a$ ending with $b$.\\

We will need the following Proposition.
\begin{Prop} \label {propmixing}
Subshift $(\Omega,X,\sigma)$ is topologically mixing.
\end{Prop}
Note that although our original map is clearly topologically mixing elements of the Markov partition are Cantor sets,
so the statement is not obvious. \\

Proof of Proposition \ref{propmixing}. By construction any Cantor set $C$ which coincides  with some element of the tower is mapped by some iterate of $F$ onto a full width substrip of some Cantor set $C_i$. As the image of any $C_i$
(including $C_1$) contains a full with substrip of $C_1$,
we get that all consecutive images of $C_i$ have Markov intersections with $C_1$.\\

It remains to prove that for any element $\Delta$ of the tower 
there is  an $n(\Delta)$ such that 
$F^n C_1$ has Markov intersection with $\Delta$ for $n >
n(\Delta)$ .  By construction any $\Delta = F^{k(\Delta)}P$ where $P$ is a full height Cantor subrectangle of some $C_i$.
So it is enough to prove that $F^n C_1$ intersects $C_i$ for $n>n(i)$. But $C_1$ contains  
a full height Cantor subset $C_{11\ldots 1 i}$ 
 and all  images of $C_1$ have Markov  intersections with   $C_{11 \ldots i}$   . That proves Proposition \ref{propmixing}.\\

Next we consider the first return map $T_{1}$ induced  by $T$ on $C_{1}$. 
Consider  the Markov partition ${ \cal MP}_1$ of $C_{1}$
\beq  \label {MPnew} 
 { \cal MP}_1 = \{C_{1i_1 \ldots i_{m}1}\}
\eeq
generated on $C_1$ by $T_1$.
By construction $T_{1}$ maps its domains (which are full height Cantor sets) onto full width substrips of $C_1$.
 Therefore the transition matrix
corresponding to the map $T_{1}$ on ${ \cal MP}_1$  consists of all $1$'s.\\

\textit{\textbf{Remark:} $[1]$ will correspond to our choice of state $[a]$ in the  symbolic dynamics in later sections.}

\ee


\section{Thermodynamic formalism,  Reduction of Theorems \ref{main} and \ref{CLT} to results for functions defined on one-sided sequences. }
\be

\item \textbf{Reduction Arguments.}

By first reducing to functions defined on one-sided sequences, we will show that our transfer operator (\ref{Transfer1}) has the spectral gap property (\ref{spectgap}) on a particular Banach Space. This property implies exponential decay of correlations and  Central Limit Theorem for functions defined on one-sided sequences. Then we can extend these results to  certain functions defined on two-sided sequences, Theorems \ref{main} and \ref{CLT}, respectively. This exchange between the two settings is a consequence of the reduction from two-sided shifts to one-sided shifts following from the classical arguments of Ruelle and Bowen,
see
\cite{Ruelle-1}, \cite{Bowen}. In the  case of an infinite alphabet, detailed reduction arguments can be found in 
  sections 4 and 5  of \cite{Young1}. \\

If we restrict our consideration to  H\"older functions on $Q$, then we are left with the proof of the spectral gap property (SGP) of the transfer operator 
 acting on a suitable space $\cal L$ of functions defined on one-sided sequences of the alphabet $\Omega$.
We do this in the next sections, following  \cite{Cyr} and \cite{Sarig-4}.

\item \textbf{Thermodynamic Formalism.}

Now by following \cite{Ruelle-1}, \cite{Bowen} we develop thermodynamic formalism on the space of
 one-sided sequences for the function $\Phi(x,y) = -log|D^{u}F|$.\\
 
Let ${\cal E}_{i} \subset E_i$ be an element of the Markov partition ${\cal MP}$. For each ${\cal E}_i$ we fix some unstable manifold $W^{u}_{0}$, and to any $(x,y) \in {\cal E}_{i}$ we let correspond $(x,y_{0}) = W^{s}(x,y) \cap W^{u}_{0}$. We define

$$
u(x,y) = \sum_{k=0}^{\infty} \Phi(F^{k}(x,y))-\Phi(F^{k}(x,y_{0})).
$$

\noindent Then we can construct a H\"older function on one-sided sequences cohomologous to $\Phi(x,y)$ in the following way,

\beq \label{POTENTIAL}
\phi(x)=\Phi(x,y) + u(F(x,y))-u(x,y).
\eeq

\noindent We call $\phi$ the potential.\\

The transfer operator, $L_{\phi}$  is defined as 

\beq \label{Transfer1}
(L_{\phi}f)(x) =\sum_{F(y)=x} e^{\phi(y)}f(y).
\eeq

In the next several sections we consider the space of functions on one-sided sequences for which we can prove the spectral gap property for
$L_{\phi}$. We denote by $X$ the space of one-sided sequences. From this point forward, points $x,y$ will be one-sided sequences belonging to $X$.

\item \textbf{Induced system}. \\ 
Just as in general case, we get SGP as a result of a particular induction procedure, see
  \cite{Cyr}, \cite{Sarig-4}, and  references
to earlier works  in \cite{Sarig-4}. In our setting we induce on $C_{1}$; here our $[a]$ is $[1]$.

The induced system on $[1]$ is $F_{1}:X_{1} \to X_{1}$ where 
$$
X_{1} := \{ x \in X\> \> | \> \> x_{0}=1, x_{i}=1 \> \> \textnormal{infinitely often} \}
$$

\noindent and

$$F_{1}(x) := F^{\varphi_{1}}(x),
$$

\noindent for 

$$\varphi_{1}(x):=1_{[1]}(x)\min\{n\geq 1: \> x_{n}=1\}.
$$

\noindent The resulting transformation can be given the structure of a Markov Shift as follows. Let

$$\overline{S}:= \{ [1,\xi_{2},...,\xi_{n-2},1 ] : 2 \leq i \leq n-3, \> \> \xi_{i} \neq 1 \} 
$$

\noindent and let $F^{\varphi_{1}}:\overline{X}\to \overline{X}$ denote the left shift on $\overline{X}=(\overline{S})^{\N}$. Then $F^{\varphi_{1}}$ is topologically conjugate to $F_{1}$. The conjugacy $\overline{\pi}:\overline{X} \to X_{1}$ is given by

$$
\overline{\pi}([1, \underline{\xi}^{0},1], [1, \underline{\xi}^{1},1], ...):= (1, \underline{\xi}^{0}, 1, \underline{\xi}^{1},...).
$$

Let
$$
\overline{\phi}:=\left ( \sum_{i=0}^{\varphi_{1}-1} \phi \circ F^{i} \right ) \circ \overline{\pi}.
$$
\noindent We call $\overline{\phi}: \overline{X} \to \R$  the induced potential.\\

\item \textbf{Gurevich Pressure.}\\

We introduce a few preliminary definitions and results.\\

\noindent Let $\phi_n(x)= \sum_{k=0}^{n-1}\phi \circ F^k(x)$ where $x=(.x_{0}x_{1} \dots)$.
Let
\beq  \label {statsum}
Z_n(\phi) = \sum_{T^n x=x; \> x_0 = 1} e^{\phi_n(x)}. 
\eeq
Then the limit,  called the {\it Gurevich Pressure},
\beq  \label {Potphi}
P_{G}(\phi) = \lim_{n \to \infty}\frac{1}{n}\log Z_n(\phi) 
\eeq
exists, see \cite{Sarig-1}, and we can calculate it explicitly in our setting. \\

\noindent Let $T_{1}$ be the first return map to $C_{1}$. Each periodic orbit of $T_{1}$ is contained in some admissible cylinder of $F$, also periodic but which has, in general, a larger period. Moreover 
there are $F$-strings of arbitrary large periods which correspond to a given $T_{1}$ period. Admissible cylinders of the same $T_{1}$ periods but with different $T_{1}$ labels do not intersect.\\

\noindent Proposition \ref{unifdist1} and uniformly bounded distortions of $D^u F$ imply that the contribution to $P_{G}(\phi)$ from each
 periodic $T_{1}$ orbit
differs from the length of the horizontal cross section of the respective two-dimensional  cylinder by a uniformly bounded factor. That implies that  the quantities $Z_n(\phi)$ are uniformly bounded from above.\\

\noindent As the measure of the Cantor set $C_{1}$ is positive, we get that $Z_n(\phi)$ are uniformly bounded from below. Thus,

$$ P_{G}(\phi)=0.$$

\noindent \textit{\textbf{Remark:} This calculation works for any $C_{i}$, so as in the general case (see \cite {Sarig-4}) the Gurevich pressure is independent of the choice of partition set $[a]$.}

\item \textbf{Spectral gap.}\\

\noindent Following  \cite {Cyr}, \cite{Sarig-4} we would like to show that $L_{\phi}$ has spectral gap on the appropriate Banach space $\mathcal{L}$, defined below in (\ref{Ldef}).  The implications 
of such a property are as follows.\\

\noindent If $L_{\phi}$ has spectral gap then it can be written as

$$
L_{\phi} = \lambda P + N
$$

\noindent where

$$
\lambda = e^{P_{G}(\phi)}, \> \> PN=NP=0, \> \> P^{2}=P
$$

\noindent and the spectral radius, $\rho$, of $N$ is less than $\lambda$. Since $\rho < \lambda$,

\beq \label{spectgap}
\dnorm{\lambda^{-n}L_{\phi}^{n}-P}_{\cal L}=\lambda^{-n} \dnorm{N^{n}}_{\cal L} \to 0
\eeq

\noindent exponentially fast as $n \to \infty$.\\
 
\noindent In our setting,

$$
\lambda = e^{P_{G}(\phi)}=1
$$

\noindent and

\begin{equation} \label{eigenmeasure}
Pf = h \int f \> d\nu
\end{equation}
where $h$ is the eigenfunction of $L_{\phi}$  and $\nu$ is the eigenmeasure of $L^{*}_{\phi}$.\\

Following \cite {Cyr}, \cite{Sarig-4} we introduce the 
 $a$-discriminant 
$$\Delta_{a}[\phi]:=\sup_{p \in \R} \{P_{G}(\overline{\phi + p}) \> \> | \> \> P_{G}(\overline{\phi + p}) < \infty \}.
$$
\noindent The Discriminant Theorem 6.7, from \cite{Sarig-4}, gives necessary and sufficient conditions for the spectral gap 
based on certain properties of the $a$-discriminant. Specifically, it links the strict positivity of the discriminant to the spectral gap property. It  involves the Gurevich pressure evaluated with respect to the induced system.
We state one of the properties relevant to our setting.

\begin{Prop} \label{discriminant}
Let $X$ be a topologically mixing TMS and suppose $\phi: X \to \R$ is a weakly H\"older continuous function such that $P_{G}(\phi) < \infty$. If for some state $a$, 
\beq \label {discrcondition}
\Delta_{a}[\phi] > 0,
\eeq
then $\phi$ has the SGP on the Banach space ${\mathcal L} $ defined below in (\ref{Ldef}).
\end{Prop}

\item The following results prove $\Delta_{1}[\phi] > 0$ by showing that for sufficiently small $p$,

$$
0 <  P_{G}(\overline{\phi + p}) < \infty.
$$

\begin{Prop} \label{DiscrPos1}
For sufficiently small $p>0$, $ P_{G}(\overline{\phi + p}) < \infty$.
\end{Prop}

Recall that

\beq  
P_{G}(\overline{\phi+p}) = \lim_{n \to \infty}\frac{1}{n}\log Z_n(\overline{\phi + p}). 
\eeq

We begin by calculating $Z_{1}(\overline{\phi + p})$:

$$
Z_{1}(\overline{\phi + p}) = \sum_{ \{x:\> Tx=x; \>\> x_{0}=1 \} } e^{\overline{\phi +p}} = \sum_{n=1}^{\infty} \left [ \sum_{ \{x:\> Tx=F^{n}x=x; \>\> x_{0}=1 \}}   e^{\overline{\phi}}e^{pn} \right ]
$$

$$
=\sum_{n=1}^{\infty} e^{pn} \left [ \sum_{ \{x:\> Tx=F^{n}x=x; \>\> x_{0}=1 \}}   e^{\overline{\phi}} \right ]
$$

$$
\leq \sum_{n=1}^{\infty} e^{pn} C_{\phi} \beta^{n} = C_{\phi}M_{1}<\infty
$$

for $p < \log\left ( \frac{1}{\beta}\right).$ \\

Note that this $\beta$ comes from the estimate in Proposition \ref{FR}.\\

Next we calculate $Z_{2}(\overline{\phi + p})$,

$$
Z_{2}(\overline{\phi + p}) = \sum_{ \{x:\> T^{2}x=x; \>\> x_{0}=1 \} } e^{\overline{\phi +p}} 
$$

$$
= \sum_{n_{1}=1}^{\infty} \sum_{n_{2}=1}^{\infty} \sum_{ \{x:\> Tx=F^{n_{1}+n_{2}}x=x; \>\> x_{0}=1 \}}   e^{\overline{\phi}}e^{p(n_{1}+n_{2})} 
$$

$$
= \sum_{n_{1}=1}^{\infty} e^{pn_{1}} \sum_{n_{2}=1}^{\infty} e^{pn_{2}}\sum_{ \{x:\> Tx=F^{n_{1}+n_{2}}x=x; \>\> x_{0}=1 \}}   e^{\overline{\phi}} 
$$

$$
\leq \sum_{n_{1}=1}^{\infty} e^{pn_{1}} \sum_{n_{2}=1}^{\infty} e^{pn_{2}}\>(C_{\phi}^{2} \beta^{n_{1}+n_{2}})
$$

$$
=\sum_{n_{1}=1}^{\infty} e^{pn_{1}} C_{\phi} \beta^{n_{1}}\sum_{n_{2}=1}^{\infty} e^{pn_{2}}\>C_{\phi} \beta^{n_{2}} \leq C_{\phi}^{2}M_{1}^{2}< \infty.
$$

Similarly we prove by induction for all $n$,

$$
Z_{n}(\overline{\phi+p}) \le C_{\phi}^{n}M_{1}^{n},
$$

and thus,

$$
P_{G}(\overline{\phi+p}) = \lim_{n \to \infty}\frac{1}{n}\log Z_n(\overline{\phi + p}) \leq \lim_{n \to \infty}\frac{1}{n} \log ( C_{\phi}^{n}M_{1}^{n})
$$

$$
=C_{\phi} + M_{1} < \infty.
$$

\begin{Prop} \label{DiscrPos2}
For $p>0$ as in Proposition \ref{DiscrPos1}, $ P_{G}(\overline{\phi + p}) >0.$
\end{Prop}

\noindent We combine the following properties from \cite{Sarig-4} with the fact that in our setting $P_{G}(\phi)=0$ for $\phi(x)$ cohomologous to $\Phi(x,y)$.

\begin{itemize}
\item[i.] $\overline{\phi + p} \geq \overline{\phi} +p$, $\forall p \in \R^{+}$
\item[ii.] If  $\phi \leq \psi$, then $P_{G}(\phi) \leq P_{G}(\psi)$
\item[iii.] $P_{G}(\phi + p) = P_{G}(\phi) + p$, $\forall p \in \R$
\item[iv.] $P_{G}(\phi)=0 \iff P_{G}(\overline{\phi})=0$
\end{itemize}

\noindent For our $p>0$, 

$$
 P_{G}(\overline{\phi + p}) \geq P_{G}(\overline{\phi} + p) = P_{G}(\overline{\phi}) + p  =p>0.
$$
\item \textbf{Banach Space.}

For all $x,y \in X$, let

$$
t(x,y) = min\{n:x_{n} \neq y_{n}\}
$$

$$
s_{1}(x,y) = \#\{0 \leq i \leq t(x,y) -1: x_{i}=y_{i}=1\}.
$$

Let $[1]$ be the collection of one-sided sequences, $x=.x_{0}x_{1}...$, such that $x_{0}=1$ .\\

As proved in  \cite{Cyr}  there is a positive function $h_0: {\bf Z}_+ \rar \bf R$ with the following properties. \\

Consider the set  of continuous functions $\{f: \parallel f \parallel_{\mathcal{L}} \> < \infty \}$ where

\bqn \label{Ldef}
\parallel f \parallel_{\mathcal{L}} \> = \sup_{b \in {\bf Z}_+} 
\frac{1}{h_0(b)}  
 \left [  \sup_{x \in [b]}   | f (x)| \> + 
 \sup_{(x,y)\in [b]; \> x \neq y }  \frac{\mid f(x) - f(y) \mid}{\theta^{s_{1}(x,y)} } \right  ].  
\eqn

 Then  $ \mathcal{L}$ is an $L_{\phi}$-invariant Banach space, and 
$L_{\phi} $ on $ \mathcal{L}$  is a bounded operator with spectral gap. Additionally the eigenfunction
$h$ of Ruelle operator  belongs to  
 $ \mathcal{L}$ and for any bounded H\"older function $\psi$ it holds that
\beq \label{important}
\psi h \in   \mathcal{L}
\eeq
Note that bounded H\"older functions belong to  $ \mathcal{L}$.

\item 

It follows from Propositions \ref{DiscrPos1} and \ref{DiscrPos2} that the discriminant is strictly positive, and thus, by Proposition \ref{discriminant}, $L_{\phi}$ has spectral gap on the Banach Space $\mathcal{L}$. As in \cite {Cyr}, \cite{Sarig-4} this implies that $(\sigma, \mu_{\phi})$ has exponential decay of correlations.\\

\begin{theorem} \label{onesided}
For $\sigma$ a one-sided full shift,  consider $\phi$ the potential defined in (\ref{POTENTIAL}),
and let $ \mu_{\phi}$  be the respective invariant measure.
 Then $(\sigma, \mu_{\phi})$ has exponential decay of correlations  for  bounded H\"older functions $f $ and $g \in L^{\infty}(\mu_{\phi})$.
Namely there exists $0<\eta_{1}<1$ such that
\beq \label{decay} 
\left | \int f \>  (g \circ \sigma^n) \> d\mu_{\phi} - \int f \> d\mu_{\phi} \int g \> d\mu_{\phi} \right | < C \parallel g \parallel_{\infty} \parallel fh \parallel_{\mathcal{L}}\eta_{1}^{n}.
\eeq
\end{theorem} 
Note that $fh \in   \mathcal{L}$ because of \ref{important}.
As in \cite {Cyr},  \cite{Sarig-4}, the subsequent result follows. 
\begin{theorem} \label{CLTonesided} ({\bf Central Limit Theorem for one-sided  shift})\\
Let $(\sigma, \mu_{\phi})$ satisfy the assumptions of Theorem \ref{onesided} and suppose that  $f \in \mathcal{L}$. If $\int f d\mu_{\phi} = 0$ and $f$ cannot be expressed as $g-g \circ \sigma$ for $g$ continuous, then there is a positive, finite constant $\mathsf{d} = \mathsf{d}(f)$ such that for every $t \in \R$,

$$
\lim_{n \to \infty} \> \> \mu_{\phi} \left \{  x: \frac{1}{\sqrt{n}} \sum_{k=0}^{n-1}  f \circ \sigma ^{k}(x) \> < \> t  \right \}= \frac{1}{\sqrt{2 \pi \mathsf{d}^{2}}} \int_{-t}^{\infty}e^{\frac{-s^{2}}{2\mathsf{d}^{2}}} \> ds.
$$

\end{theorem}

\textit{\textbf{Remark:} It follows from Theorem $1.1$ part $d$ in \cite{Cyr}  that for $g$, bounded and H\"older continuous, $P(\phi + tg)$ is real analytic in a neighborhood of $0$.} 

\item \textbf{Exponential Decay of Correlations and Central Limit Theorem for functions of two variables.}\\

One can follow arguments of  Section 4 of \cite{Young1} to reduce the estimate of the $n$-th correlation 
for functions defined on the two-sided shift to the following estimate for functions defined on the one-sided shift,
\beq \label{mainestimate} 
\norm{ L_{\phi}^{n-2k}( L_{\phi}^{2k}(f_kh))- (\int  L_{\phi}^{2k}(f_kh) dm)h},
\eeq
where $h\>dm=d\nu$ and $\nu$ is the eigenmeasure of $L^{*}_{\phi}$ as in (\ref{eigenmeasure}). Here $2k < n$, and $f_k$ is  a piecewise constant approximation of $f$ on cylinders of length $k$. Thus $f_k$ is a H\"older function bounded by $\max|f|$. From  (\ref{important}) we get $f_kh \in {\mathcal L}$.\\

As in \cite{Young1} 
 Section 4  we get that norms $\parallel L_{\phi}^{2k}(f_kh)\parallel_{{\mathcal L} }$ are uniformly bounded by a constant which only depends on
$\max |f|$.
From (\ref{mainestimate}) we get an estimate similar to (\ref{decay}) but with a different constant and a  different 
$0< \eta < 1$. That proves Theorem \ref{main}.\\

Theorem \ref{CLT} follows from arguments in Section 5 of \cite{Young1} regarding a result referred to as Theorem [G] from \cite{gordin}. Using the spectral gap property, one obtains the following estimate,

\begin{equation} \label{cltproof}
\int \norm{L^{j}_{\phi}(f_{0}\>h)}\>dm \leq C_{0}' \parallel L^{j}_{\phi}(f_{0}\>h) \parallel_{\mathcal{L}} \>\leq C_{0} \eta_{0}^{j} \parallel f_{0}\>h \parallel_{\mathcal{L}} 
\end{equation}

\noindent where, as above, $f_{0}$ is a bounded H\"older approximation of $f$ and $0 < \eta_{0} <1$. Theorem \ref{CLT} relies on showing that the key assumption in Theorem [G] holds -  finiteness of the sum of the $L^{2}$-norms of the relevant conditional expectations. Showing that this assumptions holds reduces to showing that 
the sum of estimate (\ref{cltproof}) is bounded, and thus, we just need that $f_{0}\>h \in \mathcal{L}$. This again follows from (\ref{important}).

\item {\bf Concluding Remarks.}

The study of countable Markov partitions in the 1980's originated in particular from the work of Roy Adler \cite{Adler-79}. His work motivated the use of countable Markov partitions as a tool for studying one-dimensional dynamics with
critical points, and subsequently, Henon-like systems.\\

The first author keeps warmest memories of his visit in 1990 to IBM Thomas J Watson Research Center, when he worked within the wonderful group directed by Roy Adler.\\

The authors would like to thank Omri Sarig for useful discussions.

\end{enumerate}


\end{document}